\documentclass[12pt]{article}
\usepackage{epsfig,graphics}
\def\be{\begin{equation}}
\def\bea{\begin{eqnarray*}}
\def\ee{\end{equation}}
\def\eea{\end{eqnarray*}}
\def\ba{\begin{array}}
\def\ea{\end{array}}
\def\bi{\begin{itemize}}
\def\ei{\end{itemize}}
\newtheorem{theo}{Theorem}

\newtheorem{pro}{Proposition}

\def\ZZ{{{\rm Z}\kern-.4em{\rm Z}}}
\def\RR{{{\rm I}\kern-.2em{\rm R}}}
\def\NN{{{\rm I}\kern-.2em{\rm N}}}
\def\TT{{{\rm T}\kern-.5em{\rm T}}}
\def\CC{{{\rm I}\kern-.5em{\rm C}}}

\def\mV{\mathcal{V}}

\def\mT{\mathcal{T}}
\def\mP{\mathcal{P}}

\begin{document}
\title{A counterexample concerning the $L_2$-projector
onto linear spline spaces}
\author{Peter Oswald\\
\\
\small International University Bremen\\
\small e-mail: poswald@iu-bremen.de}

\date{}
\maketitle

\begin{abstract}
For the $L_2$-orthogonal projection $P_V$ onto 
spaces of linear splines over simplicial partitions in polyhedral domains in $\RR^d$,
$d>1$, we show that in contrast to the one-dimensional case, where 
$\|P_V\|_{L_\infty\to L_\infty} \le 3$ independently of the nature of the partition,
in higher dimensions the $L_\infty$-norm of $P_V$ cannot be bounded uniformly
with respect to the partition. This fact is folklore among specialists
in finite element methods and approximation theory but seemingly has never been
formally proved. 
\end{abstract}


\section{Introduction}\label{sec1}
Variational methods based on piecewise polynomial
approximations are a workhorse in numerical
methods for PDEs and data analysis. In particular,
least-squares methods lead to the study of the $L_2$-orthogonal
projection operator $P_V\;:\;L_2(\Omega)\to V$ onto a given spline space $V$
defined on a domain $\Omega\in \RR^d$. A question of considerable interest 
is the uniform boundedness of the $L_\infty$-norm
$$
\|P_V\|_{L_\infty\to L_\infty} := \max_{f\in L_\infty(\Omega):\;\|f\|_{L_\infty}=1}
\|P_Vf\|_{L_\infty}
$$
of $P_V$ with respect to families of spline spaces $V$. If $\Omega$ is
a bounded interval in $\RR^1$, then the question has been intensively
studied for the family of spaces of smooth splines of fixed degree $r$ over 
arbitrary partitions, where Shadrin \cite{shad} has recently established that
\be\label{1D}
\|P_V\|_{L_\infty\to L_\infty} \le C(r) <\infty, \qquad r\ge 1,
\ee
for any partition. This result was known for a long time for small values of
$r$, e.g., Ciesielski  \cite{cies} proved that for linear splines one can take $C(1)=3$,
while de Boor \cite{boor} solved the case $r\le 4$.
The estimate (\ref{1D}) plays an important role in numerical analysis
and for the investigation of orthonormal spline systems such as the
Franklin system in $L_p$-based scales of function spaces, $1\le p\le \infty$.

In higher dimensions, the study of the $L_\infty$-norm
of $P_V$ arose mostly in the context of obtaining $L_p$ error estimates
in the finite element Galerkin method \cite{desl,doug}, where sufficient
conditions on the underlying partition and nodal basis $\{\phi_i\}$ 
of a finite element space $V$ are formulated under which the norms
$\|P_V\|_{L_\infty\to L_\infty}$ are bounded by a certain finite constant.
Interestingly enough, these results suggest that such conditions on
partitions resp. finite element type are essential for obtaining uniform
bounds but formal proof of their {\it necessity} was not given.

Similarly, in the theory of multivariate splines final results
on the uniform boundedness of
$\|P_V\|_{L_\infty\to L_\infty}$ could not be localized.
Recently, Ciesielski \cite{cies1} asked about the extension of his result
for linear splines \cite{cies} to the higher-dimensional case, and the 
unanimous opinion of the audience was that in higher dimensions a similar result cannot hold. 
However, other than a vague reference to unpublished work by
A. A. Privalov, no concrete proof could be found.

It is the aim of this note to provide an elementary example of triangulations 
$\mT_J$ of a square $\Omega\subset \RR^2$ into $\mathrm{O}(J)$ triangles for which
\be\label{2D}
\|P_{V(\mT_J)}\|_{L_\infty\to L_\infty} \ge J, \qquad J\ge 1,
\ee
where $V(\mT_J)$ is the space of linear $C^0$ splines
(or finite element functions) on $\mT_J$; see Theorem \ref{prop2} below. 
As $J\to\infty$, the triangulations $\mT_J$ will not satisfy the minimum angle condition, 
which is natural since for these types of triangulations
a uniform bound can easily be established. However, they satisfy the maximum angle condition,
\cite{babu}, and do not possess vertices of high valence.
The example can easily be extended to $d\ge 3$, and implies that $C(1)=\infty$
for all $d>1$. 

\section{Notation and Result}\label{sec2}
We concentrate on $d=2$, the case $d\ge 3$ is mentioned in Section \ref{sec3}.
Let $\Omega\subset \RR^2$ be a bounded polygonal domain equipped with
a finite triangulation $\mT$ into
non-degenerate closed triangles $\Delta$ satisfying the usual regularity condition that 
two different triangles may intersect at a common vertex resp. edge only.
The set of vertices of $\mT$ is denoted by $\mV_\mT$.
Let $|E|$ denote the Lebesgue measure of a measurable set $E\subset \RR^2$.
By $0<{\underline{\alpha}}_\mT\le \bar{\alpha}_\mT<\pi$ we denote the minimal
and maximal interior angle of all triangles in $\mT$.

Let 
$V(\mT)$ denote the linear space of all continuous functions
$g$ whose restriction to any of the triangles $\Delta\in \mT$
is a linear polynomial. Any $g\in V(\mT)$ has a unique representation of the form
\be\label{gP}
g=\sum_{P\in\mV_\mT} g(P)\phi_P,
\ee
where the Courant hat functions $\phi_P\in V(\mT)$ are characterized by the conditions
$\phi_P(P)=1$, and $\phi_P(Q)=0$, where $Q\neq P$ is any of the remaining vertices
of $\mT$. Thus, $\dim V(\mT) = \#\mV_\mT$, and
$$
\mathrm{supp}\,\phi_P=\Omega_P:=\cup_{\Delta\in \mT:\;P\in\Delta} \Delta.
$$
The set $\Omega_P$ corresponds to the $1$-ring neighborhood of $P$ in $\mT$,
and we denote by $\mV_P=\{Q\in \Omega_P\cap \mV:\;Q\neq P\}$ the set of all neighboring vertices 
to $P$.

The $L_2$-orthogonal projection of a function $f\in L_2(\Omega)$ onto
$V(\mT)$ is given by the unique $g:=P_{V(\mT)}f\in V(\mT)$ such that
$$
(f-g,\phi_P) = 0\qquad \forall\; P\in\mV_\mT.
$$
Here and in the sequel, $(\cdot,\cdot)$ stands for the $L_2(\Omega)$ inner product.
Using (\ref{gP}) with unknown nodal values $x_P=g(P)$ as ansatz,
this is equivalent to the linear system
\be\label{Galerkin}
\sum_{Q\in \mV} (\phi_Q,\phi_P) x_Q = (f,\phi_P)\qquad
\forall \; P\in\mV.
\ee
A simple calculation shows that
$$
(\phi_P,\phi_P) = \frac{|\Omega_P|}{6}.
$$
Similarly, if $Q\in\mV_P$ is a neighbor of $P$ then
$$
(\phi_Q,\phi_P) = \sum_{\Delta\in \mT:\;P,Q\in\Delta}\frac{|\Delta|}{12},
$$
in all other cases we have $(\phi_Q,\phi_P) = 0$. This shows
in particular that
$$
(\phi_P,\phi_P) = \sum_{Q\neq P} (\phi_Q,\phi_P) =
(1,\phi_P)/2. 
$$
I.e., if we normalize in (\ref{Galerkin}) by $(\phi_P,\phi_P)$ then
(\ref{Galerkin}) turns into a linear system 
\be\label{Galerkin1}
Ax=b,\qquad x:=(x_Q:\;Q\in\mV_\mT)^T, \quad 
b:=(b_P= \frac{(f,\phi_P) }{(\phi_P,\phi_P) } :\;P\in\mV_\mT)^T,
\ee
where $\|b\|_\infty \le 2\|f\|_{L_\infty}$, and $A:=(a_{PQ})$ satisfies 
\be\label{coeff}
a_{PQ}=\left\{\ba{ll} 1,& Q=P,\\ \frac{(\phi_Q,\phi_P) }{(\phi_P,\phi_P) } > 0, 
& Q\in\mV_P, \\ 0,& \mbox{otherwise},\ea\right.
\ee
and
\be\label{coeff1}
1=\sum_{Q\neq P} a_{PQ}=\sum_{Q\in\mV_P} a_{PQ},\qquad \forall\;P\in\mV_\mT.
\ee
Thus, the matrix $A$ is only weakly diagonally
dominant, and not strictly diagonally dominant as in the one-dimensional 
case. Otherwise, we could estimate $\|A^{-1}\|_{\infty}$ in a trivial way,
and use the inequality
\be\label{Ainv}
\|P_{V(\mT)} \|_{L_\infty \to L_\infty}\le 2\|A^{-1}\|_{\infty} = 2\max_{\|Ay\|_\infty\le 1}\,\|y\|_{\infty},
\ee
which follows from the above, in conjunction with 
the obvious equality $\|P_{V(\mT)}f\|_{L_{\infty}}=\|x\|_{\infty}$.
Let us mention without proof that (\ref{Ainv}) implies the following partial result.
\begin{pro}
If for any two neighboring vertices $P\neq Q$ from $\mV_\mT$
we have $a_{PQ}\ge c_0>0$, then
$$
\|P_{V(\mT)} \|_{L_\infty \to L_\infty}\le (1+2c_0)c_0^{-2}.
$$
\end{pro}
For triangulations satisfying the minimum angle condition
uniformly, i.e., ${\underline{\alpha}}_\mT\ge {\underline{\alpha}}_0>0$, 
this result is applicable with a $c_0$ determined solely by $\alpha_0$,
and thus covers the bounds considered in the finite element literature
\cite{desl,doug}. The main result of this note is the following
\begin{theo}\label{prop2}
For any $J\ge 1$, there is a triangulation $\mT_J$ of a square into $8J+4$ 
triangles such that
the norm of the $L_2$-projector
$P_{V(\mT_J)}$ satisfies
$$
\|P_{V(\mT_J)} \|_{L_\infty \to L_\infty}\ge 2J.
$$
Thus, for spatial dimension $d=2$ we have
$$
C(1):=\sup_{\mT} \|P_{V(\mT)} \|_{L_\infty \to L_\infty}=\infty.
$$
\end{theo}

We conjecture that in terms of the number of triangles our result is asymptotically sharp
for bounded polygonal domains in $\RR^2$, i.e.,
\be\label{2DJ}
\sup_{\mT:\;\#{\mT}\le N} \|P_{V(\mT)}\|_{L_\infty\to L_\infty} \asymp N,\qquad N\to\infty.
\ee
That $\Omega$ is a square is not crucial. 
The examples given below can easily be modified to show $C(1)=\infty$ for simplicial
partitions in higher dimensions as well.

\section{Proof of Theorem \ref{prop2}}\label{sec3}

\begin{figure}[htb]
\begin{center}
\resizebox{3in}{3in}{\includegraphics{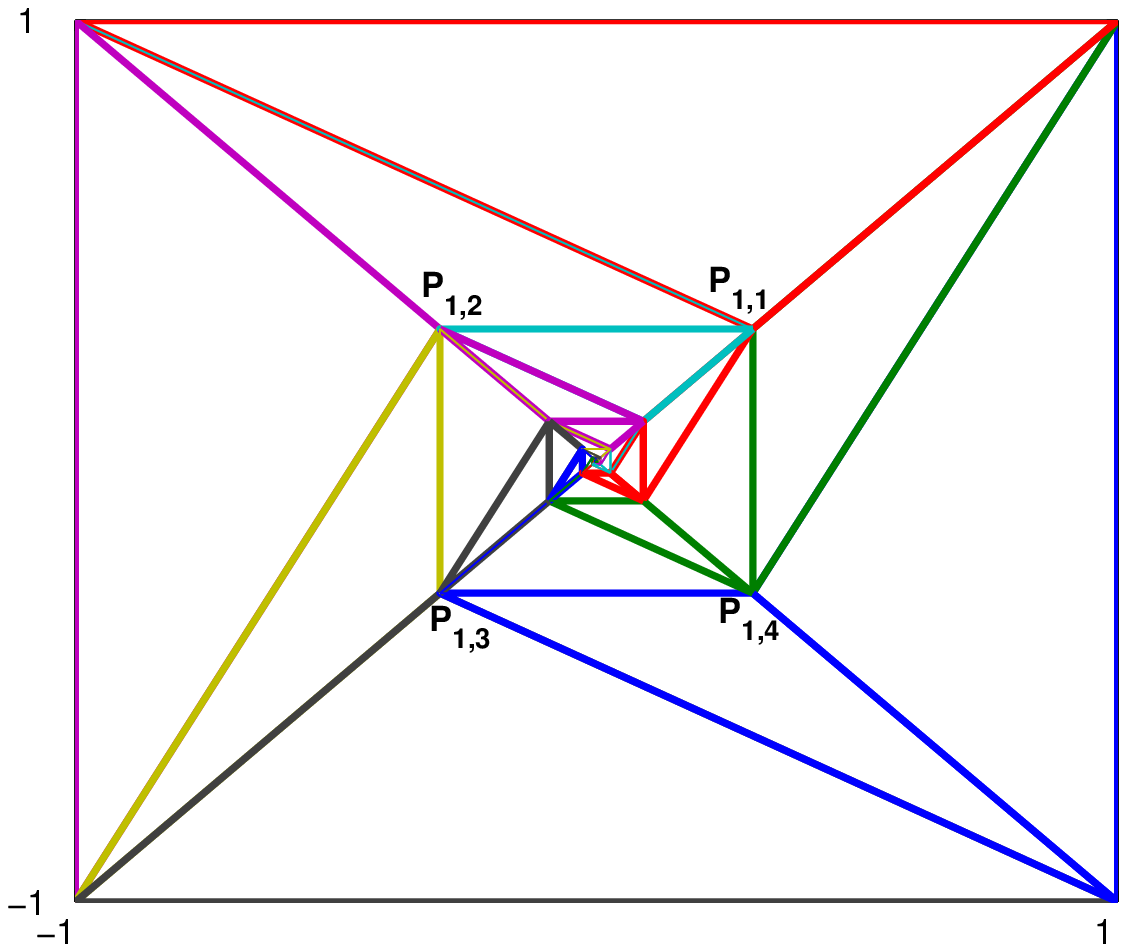}}\hspace{0.5in}
\resizebox{2in}{3in}{\includegraphics{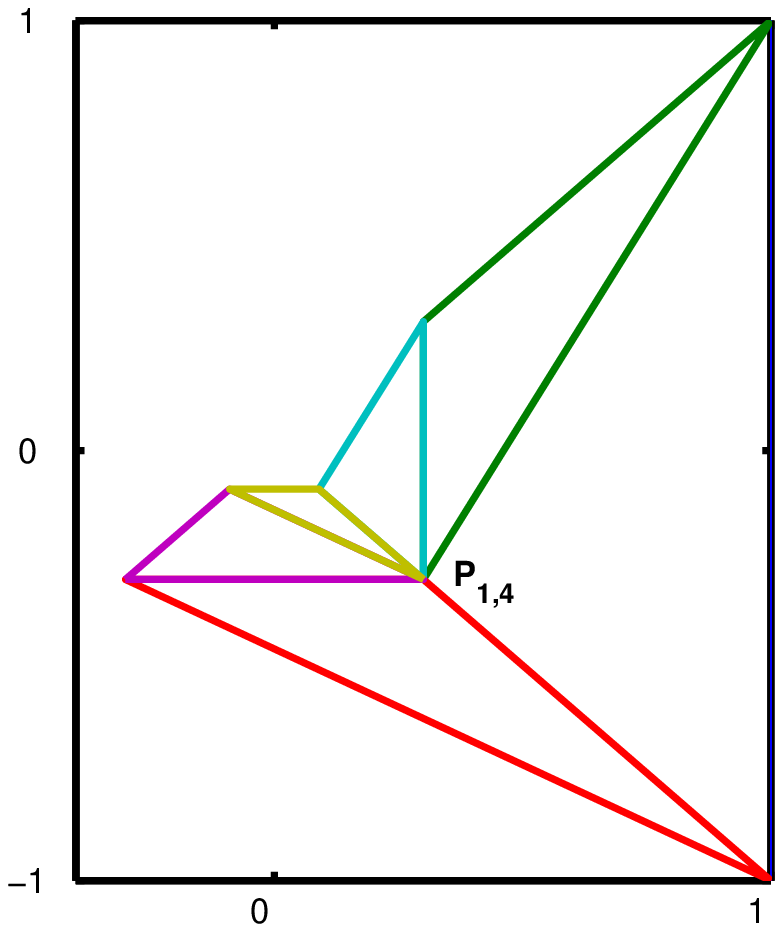}}
\end{center}
\caption{ Triangulation $\mT_J$ for $t=0.3$ (left), and typical $\Omega_{P_{j,i}}$ (right)}\label{fig1}
\end{figure}
We will use the following notation. Let
$S_a=[-a,a]^2$ be the square of side-length $2a$, $a>0$, with center
at the origin. Set $\Omega:=S_1$, with vertices
denoted (in clockwise direction) by $P_{0,i}$, $i=1,\ldots,4$,
and fix the parameter $t\in (0,1)$. The triangulation $\mT_J$, $J\ge 1$, of $\Omega$ is 
obtained by inserting the squares $S_t$, $S_{t^2}$,\ldots, $S_{t^J}$,
whose vertices will be denoted similarly by $P_{j,i}$, $i=1,\ldots,4$, $j=1,\ldots,J$,
placing an additional vertex $P_{J+1}$ at the origin, connecting
$P_{J+1}$ by straight lines with the $4$ vertices $P_{0,i}$ of $S_1$, 
and finally subdividing each of the remaining trapezoidal regions 
$P_{j-1,i-1}P_{j,i-1}P_{j,i}P_{j-1,i}$ into two triangles by connecting
$P_{j-1,i-1}$ with $P_{j,i}$ in a consistent way. The
outer rings of the resulting triangulation are shown in Fig. \ref{fig1} a). 

Let the function $f\in L_\infty(\Omega)$ with $\|f\|_{L_\infty}=1$
be defined as follows: 
$$
f(x)=(-1)^j,\qquad x\in \Omega_j := \left\{\ba{ll} S_{t^{j-1}}\backslash S_{t^j},&
j=1,\ldots, J,\\ S_{t^J},&j=J+1.\ea\right.
$$
We use the same notation as above, and consider the linear system
(\ref{Galerkin1}) corresponding to the $L_2$-orthogonal projection
$g=P_{V(\mT_J)}f$ of this $f$ onto $V(\mT_J)$. Because of uniqueness of
orthogonal projections and the rotational
symmetry of $\mT_J$ and $f$, the entries of the vector $x$ that represent 
the nodal values of $g$ corresponding
to the vertices $P_{j,i}$, $i=1,\ldots,4$, are equal, and will be denoted by
$x_j$, $j=0,\ldots,J$ (the value at the origin is denoted by $x_{J+1}$).
Moreover, the $4$ equations in (\ref{Galerkin1}) corresponding to the $4$
vertices of a square $S_{t^j}$, $j=0,\ldots,J$, can be replaced by one,
thus turning the original system $Ax=b$ of dimension $4J+5$
into a reduced tridiagonal system $\tilde{A}\tilde{x}=\tilde{b}$ of dimension
$J+2$ for the vector $\tilde{x}=(x_0,\ldots,x_{J+1})^T$.
Since the equations in the system (\ref{Galerkin1}) are invariant
under affine transformations, and because of the definition of 
$\mT_J$ via squares $S_{t^j}$ shrinking at a fixed geometric rate,
it is easy to see that, with the exception of the first and last two, 
all equations in $\tilde{A}\tilde{x}=\tilde{b}$ have the same
form:
\be\label{redA}
\alpha x_{j-1}+\beta x_j+\gamma x_{j+1}=(-1)^j\delta,\qquad j=1,\ldots,J-1.
\ee
The coefficients can be found from the triangular neighborhood $\Omega_{P_{j,i}}$ 
of any of the $P_{j,i}$, $j=1,\ldots,J-1$ (see Fig. \ref{fig1} for an illustration), and the definitions 
leading to (\ref{Galerkin1}). We do not need their exact coefficient expressions, just their limit behavior as $t\to 0$,
i.e., we will be looking for the entries of $\hat{A}:=\lim_{t\to 0} \tilde{A}$ and $\hat{b}:=\lim_{t\to 0} \tilde{b}$. 
Indeed, since for
$t\to 0$ the whole is essentially covered by the single triangle
with vertices $P_{j,i}$, $P_{j-1,i-1}$, $P_{j-1,i}$, and since $f(x)=(-1)^j$ on the latter, 
we have
$$
\alpha =1+\mathrm{O}(t),\quad \beta=1+\mathrm{O}(t),\quad \gamma =\mathrm{O}(t^2),\quad 
\delta= 2+\mathrm{O}(t).
$$
From (\ref{redA}), we obtain in the limit $t\to 0$ the equations
$$
\hat{x}_{j-1}+\hat{x}_j=2(-1)^j,\qquad j=1,\ldots,J-1,
$$
where $\hat{x}_j=\lim_{t\to\infty} x_j$ (the existence of these limits follows from the 
invertibility of the limit matrix $\hat{A}$, see below).
Similar considerations for the first and the last two
equations of $\tilde{A}\tilde{x}=\tilde{b}$ yield the remaining three equations of $\hat{A}\hat{x}=\hat{b}$
as follows:
$$
\frac32\hat{x}_{0}+\frac12 \hat{x}_1=-2,\qquad \hat{x}_{j-1}+\hat{x}_j=2(-1)^j,\quad j=J,J+1.
$$
The resulting matrix $\hat{A}$ is obviously invertible. After finding $\hat{x}_0=\hat{x}_1=-1$ from
the first two equations, forward substitution gives $\hat{x}_{j}=(2j-1)(-1)^j$, $j=2,\ldots,J+1$.
This implies that for any $\epsilon>0$ one can find a sufficiently small $t>0$ such that 
$$
\|\tilde{x}\|_\infty\ge \|\hat{x}\|_\infty-\epsilon =2J+1-\epsilon.
$$ 
This proves Theorem \ref{prop2}. Note that the above reasoning does not
work for type-I triangulations of a square obtained from a non-uniform rectangular tensor-product
partition.

We conclude with the straightforward extension of the above example to
arbitrary $d>2$. Let $e^m$ denote the $m$-th unit coordinate vector in $\RR^d$,
$m=1,\ldots,d$. As $\Omega$ we take the convex polyhedral domain with vertices
$$
P_{0,1}=e^1+e^2,\quad P_{0,2}=e^1-e^2,\quad P_{0,3}=-e^1-e^2,\quad P_{0,4}=-e^1+e^2,
$$
and $P'_m=e^m$,  $m=3,\ldots,d$. For $d=3$, this domain is a pyramid with square
base in the $xy$-plane, and tip on the $z$-axis.

A suitable simplicial partition of $\Omega$ is obtained as follows. The base
square with vertices $P_{0,i}$, $i=1,\ldots,4$, is triangulated into $\mT_J$ which depends 
on the parameters $0<t<1$ and $J$ as described
for $d=2$. The resulting triangulation (now embedded into $\RR^d$) and
its vertices are again denoted by $\mT_J$ resp. by $P_{j,i}$, $i=1,\ldots,4$, $j=0,\ldots,J$, and $P_{J+1}$. 
Then each simplex
in the associated simplicial partition $\mP_J$ of $\Omega$ is generated by the $d-2$ vertices
$P'_m$, $m=3,\ldots,d$, and the three vertices of a triangle
in $\mT_J$. The latter is called base triangle of the associated simplex. Obviously, the $d$-dimensional volume
of each simplex is proportional to the $2$-dimensional area of its base
triangle, with proportionality constant $2/d!$.
To obtain a suitable function $f\in L_\infty(\Omega)$ with $\|f\|_{L_\infty}=1$ we prescribe
values $\pm 1$ on the simplices by inheritance from the values on the
base triangles  of the above
$2$-dimensional $f$.

From the symmetry properties of $f$ and $\mP_J$, it is obvious that the $L_2$-orthogonal
projection $P_{V(\mP_J)}f$ of $f$ onto the linear spline space $V(\mP_J)$
is characterized by its value $x_{J+1}$
at the origin $P_{J+1}$, by values $x_j$ taken at the vertices $P_{j,i}$, $i=1,\ldots,4$,
where $j=0,\ldots,J$, and a common value $x'$ taken at the remaining vertices
$P'_m$, $m=3,\ldots,d$. To estimate these values which, in complete analogy to
the $2$-dimensional case, are represented by the solution vector $\tilde{x}$
of a certain linear system $\tilde{A}\tilde{x}=\tilde{b}$ (now of dimension
$J+3$), we need the limit version
$\hat{A}\hat{x}=\hat{b}$ of this linear system for $t\to 0$. We spare the reader
the elementary calculations, and state it without proof:
\bea
\hat{x}_j +\hat{x}_{j-1}+ \frac{d-2}{2} \hat{x}' &=&(-1)^j \frac{d+2}{2},\qquad j=1,\ldots,J+1,\\
\frac32\hat{x}_0 +\frac12\hat{x}_{1}+ \frac{d-2}{2} \hat{x}' &=& -\frac{d+2}{2},\\
\hat{x}_0 +\frac12\hat{x}_{1} +(1+\frac{d-3}2)\hat{x}' &=& -\frac{d+2}{2}.
\eea
From this system one easily concludes that 
$\|\hat{x}\|_\infty \ge cJ$, and consequently $\|P_{V(\mP_J)}f\|_{L_\infty}\ge cJ$
for a small enough $t>0$ which implies the desired result. 
The lower bound $cJ$ obtained does not seem to accurately
reflect the possible growth of the projector norms in $L_\infty$ as a function of
the number of simplices, for $d\ge 3$ we would rather expect an exponential rate.

\end{document}